\documentclass[11pt,leqno]{article}
\usepackage{amsmath,amssymb,amsthm,amscd}
\newcommand{\doublespace}{
   \renewcommand{\baselinestretch}{1.2}
   \large\normalsize}

\def \1{{\bf 1}}
\def \Z{\mathbb Z}
\def \C{\mathbb C}

\def \Q{\mathbb Q}

\def \wt{{\rm wt}}

\def \Res{{\rm Res}}

\def \End{{\rm End}}
\def \Hom{{\rm Hom}}

\def \<{\langle} 
\def \>{\rangle}

\def \l{\lambda }

\def \o{\omega}

\def \qed{\mbox{ $\square$}}
\def \pf {\noindent {\bf Proof:} \,}
\doublespace
\begin{document}
\newtheorem{thm}{Theorem}[section]
\newtheorem{prop}[thm]{Proposition}
\newtheorem{cor}[thm]{Corollary}
\newtheorem{lem}[thm]{Lemma}
\newtheorem{rem}[thm]{Remark}
\newtheorem{de}[thm]{Definition}
\begin{center}
{\Large {\bf The radical of a vertex operator algebra associated to a module}} \\
\vspace{0.5cm}
Chongying Dong\footnote{Supported by NSF grant 
DMS-9700923 and a research grant from the Committee on Research, UC Santa Cruz.} and Geoffrey Mason\footnote{Supported by NSF grant DMS-9700909 and a research grant from the Committee on Research, UC Santa Cruz.}\\
Department of Mathematics, University of
California, Santa Cruz, CA 95064 \\
\end{center}

\hspace{1.5 cm}
\begin{abstract} The radical of a vertex operator algebra 
associated to a module is defined and computed. 
\end{abstract}

\section{Introduction}

The study of the radical $J(V)$ for a vertex operator algebra $V$ was 
initiated in [DLMM], where we defined the radical $J(V)$ and determined 
$J(V)$ in the case $V$ is CFT type (see Section 3 for the definition
of CFT type vertex operator algebra). 
Let $M$ be an admissible $V$-module (see [DLM2] and below for
the definition). The $M$-radical $J_M(V)$ of $V$ consists of
vectors $v\in V$ such that $o(v)=0$ on $M$ where $o(v)=v_{\wt v-1}$
if $v$ is homogeneous and $o(u+v)=o(u)+o(v).$ 
In the case that $M=V,$ $J_V(V)$ is exactly $J(V)$  
which was determined to be $(L(0)+L(-1))V+J(V)_1$ in [DLMM]
where $J(V)_1$ is the weight one subspace of $J(V).$ It turns
out that a similar result is true for $J_M(V).$ We show in this paper that
$J_M(V)=(L(0)+L(-1))V+J_M(V)_{(0,1)}$ where $J_M(V)_{(0,1)}$
is the intersection of $J_M(V)$ with $V_0\oplus V_1.$ 
Although the method
for determining $J_M(V)$ is similar to that for determining $J_V(V)$
in [DLMM], the argument here is more complicated. The reason
is that $V$ has a vacuum ${\bf 1}$ but $M$ does not have in general.
We expect that the concept of  $M$-radical $J_M(V)$ of $V$ will play
a very important role in the theory of vertex operator
algebra. 

The second main result in this paper is a criterion for 
irreducibility of an admissible $V$-module $M$ (see Proposition \ref{p5.3}). 
The result says that $M$ is irreducible if and only if each
homogeneous subspace is irreducible $\hat V(0)$-module or $S_M(V)$-module
(see Section 4 for the definition of $\hat V(0)$ and $S_M(V)$).  
We also formulate this result
in terms of the theory of associative algebras $A_n(V)$ 
developed in [DLM4]. This result is
important in the study of dual pair associated to a vertex operator
algebra and an automorphism group (cf. [DLM1]). 

Both results are extended to twisted modules. In particular
we also define the radical $J_V(M)$ for an admissible
$g$-twisted $V$-module for an automorphism $g$ of $V$ of finite
order  and determine $J_V(M)$ precisely. A similar
criterion of irreducibility of $M$ is obtained too in terms 
of certain Lie algebra $S_M(V^0)$ (see Section 5) 
and associative algebra $A_{g,n}(V)$ [DLM5]. 

\section{Preliminary}
 
Let $(V,Y,{\bf 1},\omega)$ be a vertex operator algebra (see [B] and [FLM]).
We shall use commuting formal variables $z,z_0,z_1,z_2.$ We shall also
use delta-function  $\delta(z)=\sum_{n\in\Z}z^n$ whose
elementary properties can be found in [FLM]. 

First recall from [FLM], [Z], [DLM2] the definitions of weak module,
admissible module, and ordinary module for a vertex operator algebra $V.$
 A {\em weak module}
 $M$ for $V$ is a vector space equipped with a linear map
$$\begin{array}{l}
V\to (\mbox{End}\,M)[[z^{-1},z]]\label{map}\\
v\mapsto\displaystyle{Y_M(v,z)=\sum_{n\in \Z}v_nz^{-n-1}\ \ (v_n\in\mbox{End}\,M)}
\mbox{ for }v\in V\label{1/2}
\end{array}$$
satisfying the following conditions for $u,v\in V$, 
$w\in M$:
\begin{eqnarray*}\label{e2.1}
& &v_nw=0\ \ \  				
\mbox{for}\ \ \ n\in \Z \ \ \mbox{sufficiently\ large};\label{vlw0}\\
& &Y_M({\bf 1},z)=1;\label{vacuum}
\end{eqnarray*}
 \begin{equation}\label{jacobi}
\begin{array}{c}
\displaystyle{z^{-1}_0\delta\left(\frac{z_1-z_2}{z_0}\right)
Y_M(u,z_1)Y_M(v,z_2)-z^{-1}_0\delta\left(\frac{z_2-z_1}{-z_0}\right)
Y_M(v,z_2)Y_M(u,z_1)}\\
\displaystyle{=z_2^{-1}\delta\left(\frac{z_1-z_0}{z_2}\right)
Y_M(Y(u,z_0)v,z_2)}.
\end{array}
\end{equation}
Here and below  $(z_i-z_j)^n$ for $n\in \C$ is to be expanded in nonnegative
powers of the second variable $z_j.$ 

This completes the definition. We denote this weak module by
$(M,Y_M)$ (or briefly by $M$).

An {\em ordinary} $V$-{\em module} is a  weak $V$-module which
carries a $\C$-grading
$$M=\bigoplus_{\lambda \in{\C}}M_{\lambda} $$
such that $\dim M_{\l}$ is finite and $M_{\l+n}=0$
for fixed $\l$ and $n\in {\Z}$ small enough. Moreover one requires that
$M_{\l}$ is the $\l$-eigenspace for $L(0):$
$$L(0)w=\l w=(\mbox{wt}\,w)w, \ \ \ w\in M_{\l}$$
where $L(0)$ is the component operator of 
$Y_M(\omega,z)=\sum_{n\in\Z}L(n)z^{-n-2}.$

An {\em admissible} $V$-{\em module} is 
a  weak $V$-module $M$ which carries a
${\Z}_{+}$-grading 
$$M=\bigoplus_{n\in {\Z}_{+}}M(n)$$
($\Z_+$ is the set all nonnegative integers) such that if $r, m\in {\Z} ,n\in {\Z}_{+}$ and $a\in V_{r}$
then
$$a_{m}M(n)\subseteq M(r+n-m-1).$$
Note that  any ordinary 
module is an admissible module.

A vertex operator algebra $V$ is called {\em rational} if any 
admissible module is a direct sum of irreducible admissible 
modules. It was proved in
[DLM3] that if $V$ is rational then there are only
finitely many inequivalent irreducible admissible  modules and
each irreducible admissible module is an ordinary module.

The following proposition can be found in [L2] and [DM].

\begin{prop}\label{p1} Any irreducible weak $V$-module $M$ is spanned
by $\{u_nw|u\in V,n\in\Z\}$ where 
$w\in W$ is any fixed nonzero vector. 
\end{prop}

Let $M$ be a weak $V$-module. We define the $M$-radical of $V$ to be
\begin{equation}\label{3.1}
J_M(V)=\{v\in V|o(v)|_M=0\}
\end{equation}
where $o(v)=v_{\wt v-1}$ for 
homogeneous $v\in V$ and $o(u+v)=o(u)+o(v).$   If $M=V$ this
is precisely the definition of radical of $V$ given in [DLMM]. 
($J_V(V)$ was  denoted by $J(V)$ in [DLMM]). 
If $M=\oplus_{n\geq 0}M(n)$ is an admissible module then $o(v)M(n)\subset
M(n)$ for all $n\in\Z.$ 

Recall from [DLMM] that $V$ is of CFT type if $V$ is simple
and $V=\oplus_{n\geq 0}$ with $V_0$ one-dimensional.
It was proved in [DLMM] 
that if $V$ is of CFT type then 
$J_V(V)$ is equal to $(L(0)+L(-1))V+J_V(V)_1$ where 
$J_V(V)_1=V_1\cap J(V).$ Here we prove
a similar result for $J_M(V)$ for 
any admissible module $M$ with the same assumption on $V.$ 

\section{Determination of $J_M(V)$}

We need several lemmas. 

\begin{lem}\label{l4.a1} Let $V$ be a simple vertex operator algebra
and $M$ a weak $V$-module. Let $u\in V$ such that the vertex operator 
$Y_M(u,z)$ on $M$ involves only either finitely many positive powers or finitely many negative 
powers of $z$ then $u\in V_0.$
\end{lem}

\pf The proofs in the two cases are similar. We only deal with the case
that $Y_M(u,z)$ involves only finitely many positive powers of $z.$ 
We first prove that 
$$Y_M(u,z_1)Y_M(v,z_2)=Y_M(v,z_2)Y_M(u,z_1)$$
for all $v\in V.$  By (7.24) of [DL] (also see [FLM]) there exists
a nonnegative integer $n$ such that 
\begin{equation}\label{new}
(z_1-z_2)^nY_M(u,z_1)Y_M(v,z_2)=(z_1-z_2)^nY_M(v,z_2)Y_M(u,z_1).
\end{equation}
Since each factor in (\ref{new}) involves only finitely many 
positive powers of $z_1$ we multiply (\ref{new}) by $(z_1-z_2)^{-n}$
to obtain  $Y_M(u,z_1)Y_M(v,z_2)=Y_M(v,z_2)Y_M(u,z_1).$

From $[Y_M(\omega,z_1),Y_M(u,z_2)]=0$ we see that 
$$0=[L(-1),Y_M(u,z)]=Y_M(L(-1)u,z).$$ From the Jacobi identity (\ref{jacobi}) 
we have the associator formula (see Chapter 8 of [FLM]): for $a,b\in V$ and $w\in M$ there exists
a nonnegative integer $n,$ which depends on $a$ and $w$ only, such that 
$$(z_0+z_2)^nY_M(Y(a,z_0)b,z_2)w=(z_0+z_2)^nY_M(a,z_0+z_2)Y_M(b,z_2)w.$$
So if $b=L(-1)u$ then $Y_M(b,z_2)=0$ on $M$ and 
$$(z_0+z_2)^nY_M(Y(a,z_0)b,z_2)w=0$$ 
or that 
$$Y_M(Y(a,z_0)b,z_2)w=0.$$ 
This shows that $Y_M(a_mb,z)=0$ on $M$ for any $a\in V$ and $m\in\Z.$ Assume
that $b\ne 0.$ Since
$V$ is simple then the span of $a_mb$ for $a\in V$ and $m\in \Z$ 
is the whole $V$ by Proposition \ref{p1}.
As a result we have
$Y_M(v,z)=0$ for every $v\in V.$  This is a contradiction
as $Y_M({\bf 1},z)=id_M.$ Thus $b=L(-1)u=0.$ Since $L(-1): V_n\to V_{n-1}$ is
injective if $n\ne 0$ (cf. [L1] and [DLiM]) we immediately have that
$u\in  V_0,$ as required.
\qed

\begin{lem}\label{l4.a2} Let $V$ be a vertex operator algebra of CFT type. 
Let $v\in V_n$ with $n\geq 2$ such that $L(1)v=0.$ Then
$v\not\in J_M(V).$ 
\end{lem}

\pf Assume that $v\in J_M(V).$ Then $o(v)=v_{n-1}=0$ on $M.$ 
Using the relation $[L(-1),v_m]=-mv_{m-1}$ we see that
$v_k=0$ for $0\leq k\leq n-1.$ Thus for any $u\in V$ we have
$$0=[v_i,u_{-i}]=\sum_{t=0}^i\binom{i}{t}(v_tu)_{-t}$$
for $i=0,...,n-1.$ This shows that 
$$(v_iu)_{-i}=0$$ for $i=0,...,n-1.$
Using the relation $[L(-1),a]=-ma_{m-1}$  repeatedly for $a\in V$ 
gives 
$$(v_iu)_k=0$$
for $i=1,..,n-1$ and $k\leq -i.$
Thus $Y_M(v_iu,z)$ involves only finitely many
positive powers of $z.$ It follows from Lemma
\ref{l4.a1} that $v_iu\in V_0$ for $i=1,..n-1.$ If $u$ is homogeneous 
of weight $s\geq 0$ then the weight $v_{n-1}u$ again is $s.$
Thus if $s>0$ then $v_{n-1}u=0.$ If $s=0$ then $u$ is a multiple
of ${\bf 1}$ and again $v_{n-1}u=0.$ Thus $v\in J_V(V).$ 

On the other hand $J_V(V)=J_V(V)_1+(L(-1)+L(0))V$ (Theorem 1 of
[DLMM]). It is clear that $v\not\in J_V(V).$ This is a contradiction.
\qed 

\begin{lem}\label{l4.a3} Let $V$ be a vertex operator 
algebra of CFT type and $M$ a weak $V$-module.  
Let $v\in V_1$ such that  $o(v)$ is a constant on $M.$
Then $v\in J_V(V).$ Moreover if $M$ is irreducible then $o(v)$ is a constant 
on $M$ if and only if $v\in J_V(V).$
\end{lem}

\pf For any $u\in V$ and $n\in \Z$ we have 
$$0=[o(v),Y_M(u,z)]=[v_0,Y_M(u,z)]=Y_M(v_0u,z)$$
As in the proof of Lemma \ref{l4.a1} we conclude that $v_0u=0$ for
all $u\in V.$ That is, $v\in J_V(V).$

If $v\in J_V(V)$ then again we have $[o(v),Y_M(u,z)]=Y_M(v_0u,z)=0.$ 
If $M$ is irreducible then $M$ has countable dimension. Let $\Hom_V(M,M)$
denote the set of all $V$-homomorphism from $M$ to itself. Then $\Hom_V(M,M)$
is a division ring over $\C.$ Let $w\in M$ be any nonzero vector.
Then $f\mapsto f(w)$ gives a bijection from $\Hom_V(M,M)$ to $\Hom_V(M,M)w$
which has countable dimension. Thus $\Hom_V(M,M)$ has countable dimension.
Since any division ring over $\C$ with countable dimension is $\C$ itself
(cf. [DLM3]) we conclude that $\Hom_V(M,M)=\C.$ 
we now see immediately that $o(v)$ is a constant on $M.$
\qed

We can now determine the radical $J_M(V)$ precisely.
\begin{thm}\label{t3.1}
Suppose that $V$ is a vertex operator algebra of CFT type.
Then for any  admissible $V$-module $M$ we have 
$$J_M(V)=(L(0)+L(-1))V+J_M(V)_{(0,1)}$$
where $J_M(V)_{(0,1)}=(V_0+V_1)\cap J_M(V).$ Moreover, if
$a=a^0+a^1\in  J_M(V)_{(0,1)}$ with $a^i\in V_i$ then
$a^1\in J_V(V).$ That is, the image of the projection of
$J_M(V)_{(0,1)}$ into $V_1$ is contained in $J_V(V).$ 
\end{thm}

\pf The proof of this theorem is similar to that of Theorem 1 of
[DLMM]. The conclusion $(L(0)+L(-1))V+J_M(V)_{(0,1)}\subset J_M(V)$ is clear.

First we recall a result from [DLiM] (Corollary 3.2). As a module
for $sl(2,\C)=\<L(-1),L(0),L(1)\>,$ $V$ is a direct sum of highest weight
modules $X(\mu)$ with highest weights $\mu>$ ($\mu>0$), 
the trivial module and the projective cover $P(1)$ of $X(1).$ 
Thus for any $x\in J_M(V)$ we can write 
$$x=\sum_{n=0}^mL(-1)^nu^n$$
where each $u^n$ either is in $V_1$ or satisfies $L(1)u^n=0.$ 
We assume that $u^m\ne 0.$ We prove by induction
on $m$ that $x$ lies in  $J_M(V)_{(0,1)}+(L(0)+L(-1))V.$

Suppose first that $m=0.$ Then $x=u^0.$ Write $x=\sum_{i\geq 0}x^i$
where $x^i\in V_i$ and $L(1)x^i=0$ if $i\ne 1.$ 
Since $o(x)=0$ on
$M$ we have
$$0=[L(1),o(x)]=\sum_{i\geq 0}[L(1),o(x^i)]=\sum_{i\geq 2}2(i-(i-1)/2-1)x^i_i$$
on $M$ 
where we have used the fact that $[L(1),o(x^i)]=0$ for $i\leq 1$
(see Lemma 2.5 of [DLMM]). That is,
$$\sum_{i\geq 2}(i-1)x^i_i=0.$$
Thus
$$0=\sum_{i\geq 2}(i-1)[L(1),[L(-1),x_i^i]]=\sum_{i\geq 2}(i-1)[L(1),-ix^i_{i-1}]=-\sum_{i\geq 2}(i-1)^2ix_i^i$$
on $M.$ Continuing in this way we get
$$\sum_{i\geq 2}(i-1)^ki^{k-1}x^i_i=0$$
for all $k\geq 1.$ It follows that each $x^i_i=0$ for all
$i\geq 2.$ Using the relation $[L(-1),u_n]=-nu_{n-1}$ shows inductively 
that $x^i_j=0$ for $j=0,...,i.$ Thus $x^i\in J_M(V).$ If $x^i\ne 0$ then
$x^i$ is not in $J_M(V)$ by Lemma \ref{l4.a2}.
This is a contradiction. Thus $x^i=0$ for all $i\geq 2.$ 

So $x=x^0+x^1.$ Since $x^0$ is a multiple of ${\bf 1}$ we see that 
$o(x^1)=-o(x^0)$ is a constant on $M.$ By Lemma \ref{l4.a3}
$x^1\in J_V(V).$ This proves the result for $m=0.$ 

For $m>0$ set $a=L(-1)^{m-1}u^m$ and $b=\sum_{n=0}^{m-1}L(-1)^nu^n.$ Thus 
$x=L(-1)a+b.$ From $(L(0)+L(-1))a\in J_M(V)$ we have
$$0=o(x)=o(L(-1)a)+o(b)=-o(L(0)a)+o(b)=o(b-L(0)a).$$
Note that $L(0)a=(m-1)L(-1)^{m-1}u^m+L(-1)^{m-1}L(0)u^m$ so that
$$b-L(0)a=\sum_{n=0}^{m-2}L(-1)^nu^n-L(-1)^{m-1}((m-1)u^m+L(0)u^m-u^{m-1})$$
lies in $J_M(V).$ Since either $L(0)u^m\in V_1$ or
$L(1)L(0)u^m=0,$  we conclude by induction that $b-L(0)a$ lies in 
$J_M(V)_{(0,1)}+(L(0)+L(-1))V.$ But then the same is true for $x=b
-L(0)a+(L(0)+L(-1))a.$ This completes the proof of the theorem. 
\qed

\section{A criterion for irreducibility}

In this section we give a criterion for irreducibility of an admissible module
for an arbitrary vertex operator algebra $V$ which we do not assume to be 
simple. 
We consider the quotient space
$$\hat V={\C}[t,t^{-1}]\otimes V/D{\C}[t,t^{-1}]\otimes V$$
where $D=\frac{d}{dt}\otimes 1+1\otimes L(-1).$ Denote by $v(n)$ the image
of $t^n\otimes v$ in $\hat V$ for $v\in V$ and $n\in \Z.$ Then $\hat V$
is $\Z$-graded by defining the degree of $v(n)$ to be $\wt v-n-1$ if $v$
is homogeneous. Denote the homogeneous subspace of degree $n$ by $\hat V(n).$
The space $\hat V$ is, in fact, a $\Z$-graded Lie algebra with bracket
$$[a(m), b(n)]=\sum_{i=0}^{\infty}\binom{m}{i}a_ib(m+n-i)$$
for $a,b\in V$ 
(see [B], [L2] and [DLM3]). Note that $\hat V(0)$ is a subalgebra of
$\hat V$ and is isomorphic to $V/(L(-1)+L(0))V$ whose Lie
bracket is given by
$$[a,b]=\sum_{n=0}^{\wt a-1}\binom{\wt a-1}{n}a_nb$$
for homogeneous $a,b\in V.$

Let $M$ be an admissible $V$-module. Then the map from $\hat V$ to
$\End M$ by sending $v(m)$ to $v_m$ is a Lie algebra homomorphism
(cf. [L2] and [DLM3]). In particular, the restriction of this map to $\hat V(0)$
gives a Lie algebra homomorphism from $ V/(L(-1)+L(0))V$ to
$\End M.$ The kernel  of this
map is exactly the $M$-radical $J_M(V).$ Set
$$S_M(V)=V/J_M(V).$$
Then $S_M(V)$ is a quotient Lie algebra of $V/(L(-1)+L(0))V$ by Theorem
\ref{t3.1} and acts on $M$ faithfully.  

\begin{lem}\label{la1} Let $V$ be a finite dimensional vertex operator
algebra. Then $V=V_0$ is a commutative associative algebra such
that $Y(a,z)b=ab$ for $a,b\in V.$
\end{lem}

\pf Since $L(-1)$ is injective on $\sum_{n>0}V_n$ (see [L1] and
[DLiM]) we observe
that $\sum_{n>0}V_n=0.$ In particular, $\omega=0$ and $L(0)=0.$  
This shows that $V=V_0.$ 

It is clear now that $a_n=0$ for $a\in V_0$ and $n\ne -1.$ This
implies that $Y(a,z)b=a_{-1}b.$ The reader can verify that
$ab=a_{-1}b$ defines a commutative associative algebra structure
on $V_0$ (see [B] and [L2]). 
\qed 

\begin{lem}\label{l1} Let $V$ be a vertex operator algebra and
$M=\oplus_{n\geq 0}M(n)$ an admissible $V$-module
with $M(0)\ne 0.$ Then $M$ is not equal to $\oplus_{n=0}^kM(n)$ for any
$k\geq 0$ unless $V$ is finite dimensional. 
\end{lem}

\pf If  $M=\oplus_{n=0}^kM(n)$ such that $M(k)\ne 0.$ Take a nonzero 
$u\in M(k).$ Then from the definition of admissible module 
$L(-1)u\in M(k+1)=0.$ Thus $u$ is a vacuum-like vector and the submodule
$W$ of $M$ generated by $u$ is isomorphic to the adjoint module $V$ [L1].
Since $L(-2)u\in M(n+2)=0$ we see that $\omega=0$ and $V=V_0$ is finite
dimensional. \qed

Now we use Proposition \ref{p1} and Lemma \ref{l1} to give a criterion
for irreducibility of an admissible module. 
\begin{prop}\label{p5.3} Let $V$ be a vertex operator algebra
with $\omega\ne 0.$ An admissible $V$-module $M=\oplus_{n\geq 0}M(n)$ with $M(0)\ne 0$
is irreducible if and only if each $M(n)$ is an irreducible $\hat V(0)$-module,
or each $M(n)$ is an irreducible $S_M(V)$-module.
\end{prop}

\pf We have already mentioned that $M$ is a module for $\hat V$
under the action $v(n)\mapsto v_n$ where $Y_M(v,z)=\sum_{n\in\Z}v_nz^{-n-1}$
are vertex operators on $M$ for $v\in V.$ First we assume that $M$ is
irreducible. By Proposition \ref{p1}
$M=\hat V\cdot w$ for any nonzero vector $w$ of $M$. Now take
$w\in M(n).$ Then $M(k)=\hat V(k-n)\cdot w.$ In particular,
$M(n)=\hat V(0)\cdot w.$ Thus $M(n)$ is an irreducible $\hat V(0)$-module.

Conversely suppose each $M(n)$ is an irreducible $\hat V(0)$-module. 
From the proof of Lemma 1.2.1 of [Z] we see that $L(0)$ acts on
each $M(n)$ as a scalar. Let $W$ be any nonzero submodule of $M.$ Then
$$W=\oplus_{n\geq 0}W(n)$$
where $W(n)=M(n)\cap W.$ From Lemma \ref{l1} and the injectivity of
$L(-1)$ on $M(n)$ for all large $n$ (cf. [L1] and [DLiM]) 
we see that $W(n)\ne 0$ for all
large $n.$ Note that each $W(n)$ is a submodule of $M(n)$
for $\hat V(0).$ So $W(n)=M(n)$ for all large $n$ as $M(n)$
is irreducible $\hat V(0)$-module. If $W\ne M$ then the quotient
$M/W$ is an admissible $V$-module with only finitely many homogeneous
subspaces. This is a contradiction by Lemma \ref{l1} unless $V$ is
finite dimensional. Thus $W=M$ and $M$ is irreducible. 
\qed

\begin{rem} {\rm In the case $V=V_0$ is a commutative associative
algebra, the assertion in Proposition \ref{p5.3} is false. For example,
if we take $V=\C,$ then $M=\sum_{n\geq 0}M(n)$ is a $V$-module
with each $M(n)=V.$ Clearly, each $M(n)$ is an irreducible $S_M(V)$-module
but $M$ is not irreducible under $V.$ }
\end{rem} 

The result discussed Proposition \ref{p5.3} can be also formulated in terms 
of the theory of associative algebra $A_n(V)$ developed in [DLM4]. 

Let $O_n(V)$ be the linear span of all $u\circ_n v$ and $L(-1)u+L(0)u$
where for homogeneous $u\in V$ and $v\in V,$
$$u\circ_n v=\Res_{z}Y(u,z)v\frac{(1+z)^{\wt u+n}}{z^{2n+2}}.$$
Define the linear space $A_n(V)$ to be the quotient $V/O_{n}(V).$ 
We also define a second product $*_n$ on $V$ for $u$ and $v$ as
above: 
$$u*_nv=\sum_{m=0}^{n}(-1)^m\binom{m+n}{n}\Res_zY(u,z)\frac{(1+z)^{\wt\,u+n}}{z^{n+m+1}}v.$$
Extend linearly to obtain a bilinear product  on $V$ which coincides with that
of Zhu [Z] if $n=0.$ The following theorem was proved in [DLM4]; In the case $n=0$
it was proved previously in [Z]. 

\begin{thm}\label{theorem:2.1} Let $M=\sum_{n\geq 0}M(n)$ be an admissible
$V$-module with $M(0)\ne 0.$ Then

{\rm (i)} The product $*_n$ induces 
an associative algebra structure on $A_n(V)$ 
with the identity $\1+O_n(V).$ Moreover $\o+O_n(V)$ is a central
element of $A_n(V).$ 

{\rm (ii)} The identity map on $V$ induces an onto algebra homomorphism
from $A_n(V)$ to $A_m(V)$ for $0\leq m\leq n.$ 

{\rm (iii)} The map $u\mapsto o(u)$ gives a representation of $A_n(V)$ on $M(i)$ for $0\leq i\leq n.$ 
Moreover, $V$is rational if and only if $A_n(V)$ are finite dimensional
semisimple algebras for all $n.$  
\end{thm}

Note that both the actions of $A_{n}(V)$ and $S_M(V)$ on 
$\sum_{0\leq m\leq n}M(m)$ are given by $v\mapsto o(v).$ 
Combining Proposition \ref{p5.3} and Theorem \ref{theorem:2.1} immediately
gives
\begin{thm}\label{t6} Assume that the Virasoro element $\omega$ of
$V$ is nonzero. Then an admissible
$V$-module $M=\oplus_{n\geq 0}M(n)$ with $M(0)\ne 0$
is irreducible if and only if each $M(n)$ is an irreducible $A_n(V)$-module.
\end{thm}

\section{Twisted case}

This section is an analogue of Section 4 for a twisted
module $M.$ We will omit a lot of details and refer the reader
to the previous sections when it is clear how the corresponding 
proofs and arguments before carry out in this case. 

First we give  definitions of various twisted modules following [FLM] and
[DLM3]. Let $g$ be an automorphism of $V$ of order $T.$ Then we have eigenspace
decomposition $V=\sum_{k=0}^{T-1}V^k$ where 
$V^k=\{v\in V|gv=e^{\frac{-2\pi i k}{T}}v\}.$  Then $V^0$  
is a vertex operator subalgebra of $V$ with the same Virasoro vector.

A {\em weak} $g$-{\em twisted} $V$-{\em module} $M$ is a vector space equipped 
with a linear map
$$\begin{array}{l}
V\to (\End\,M)\{z\}\\
v\mapsto\displaystyle{ Y_M(v,z)=\sum_{n\in\Q}v_nz^{-n-1}\ \ \ (v_n\in
\End\,M)}
\end{array}$$
such that for all $0\leq r\leq T-1,$ $u\in V^r$, $v\in V,$ 
$w\in M$,
\begin{eqnarray*}
& &Y_M(u,z)=\sum_{n\in \frac{r}{T}+\Z}u_nz^{-n-1} \\ 
& &u_lw=0\ \ \  				
\mbox{for}\ \ \ l>>0\\
& &Y_M({\bf 1},z)=id_M;
\end{eqnarray*}
 $$
\begin{array}{c}
\displaystyle{z^{-1}_0\delta\left(\frac{z_1-z_2}{z_0}\right)
Y_M(u,z_1)Y_M(v,z_2)-z^{-1}_0\delta\left(\frac{z_2-z_1}{-z_0}\right)
Y_M(v,z_2)Y_M(u,z_1)}\\
\displaystyle{=z_2^{-1}\left(\frac{z_1-z_0}{z_2}\right)^{-r/T}
\delta\left(\frac{z_1-z_0}{z_2}\right)
Y_M(Y(u,z_0)v,z_2)}.
\end{array}
$$
It is clear if $g=1$ this reduces the definition of weak module
in Section 3. 

An {\em ordinary} $g$-{\em twisted $V$-module} is
a weak $g$-twisted $V$-module $M$ 
with a $\C$-grading induced by the eigenvalues of $L(0):$
$$M=\bigoplus_{\lambda \in{\C}}M_{\lambda} $$
where $M_{\l}=\{w\in M|L(0)w=\l w\},$ $\dim M_{\l}$ is finite and for fixed $\l,$ $M_{\frac{n}{T}+\l}=0$ for all small enough integers $n.$

An {\em admissible $g$-twisted $V$-module} is a  weak $\frac{1}{T}{\Z}$-graded
$g$-twisted $V$-module $M$ 
$$M=\bigoplus_{n=0}^{\infty}M(n/T)$$
such that $M(0)\ne 0$ and that
$v_{m}M(n/T)\subseteq M(n/T+\wt v-m-1)$ for homogeneous $v\in V.$
Clearly, an ordinary $g$-twisted $V$-module is an admissible
$g$-twisted $V$-module. 

\begin{rem} {\rm From the definition we see that any weak (admissible, ordinary)
$g$-twisted $V$-module is a weak (admissible, ordinary) $V^0$-module.}
\end{rem} 

Let $M$ be an admissible $g$-twisted $V$-module. 
For homogeneous $v\in V$ we
denote $o(v)=v_{\wt v-1}$ on $M$ and extend it linearly to whole $V,$ 
as before.
Then it is immediate from the definition that 
$o(v)=0$ for $v\in V^1\oplus\cdots \oplus V^{T-1}.$ Since $M$ is
an admissible $V^0$-module we consider the $M$-radical of $V^0$ given in
(\ref{3.1}). By Theorem \ref{t3.1} we have
\begin{thm}
Suppose that $V$ is a vertex operator algebra of CFT type.
Then for any  admissible $g$-twisted $V$-module $M$ we have 
$$J_M(V^0)=(L(0)+L(-1))V^0+J_M(V^0)_{(0,1)}$$
Moreover, if
$a=a^0+a^1\in  J_M(V^0)_{(0,1)}$ with $a^i\in V^0_i$ then
$a^1\in J_{V^0}(V^0).$ 
\end{thm}

Proposition \ref{p5.3} still holds in this case.
\begin{prop}\label{p6.3} Let $V$ be a simple vertex operator algebra with
$\omega\ne 0.$ 
An admissible $g$-twisted 
$V$-module $M=\oplus_{n\geq 0}M(n/T)$ with $M(0)\ne 0$
is irreducible if and only if each $M(n/T)$ is an irreducible 
$S_M(V^0)$-module.
\end{prop}

\pf If $M$ is irreducible then
one can show that the analogue of Proposition \ref{p1} is true. That is,
$M=\{u_nw|u\in V,n\in\frac{1}{T}\Z\}$ for any nonzero $w\in M.$ 
Thus for any nonzero $w\in M(n/T)$ we have $\{o(v)w|v\in V^0\}=M(n/T).$
So $M(n/T)$ is an irreducible $S_M(V^0)$-module.

Note that for each $k=0,...,T-1,$ $M^k=\oplus_{n\geq 0}M(n+k/T)$ is an admissible
$V^0$-module. If all $M(n/T)$ are irreducible
$S_M(V^0)$-modules then $M^k$ is an irreducible
admissible $V^0$-modules for $k=0,...,T-1$ by Proposition
\ref{p5.3}. Using the associativity of vertex operators on $M$ we show
that if $Y_M(v,z)w=0$ for some nonzero $v\in V$ and $w\in M$ then
$Y_M(u,z)=0$ for all $u\in V$ (cf. Proposition 11.9 of [DL]: here we use
the assumption that $V$ is simple).  Since $\sum_{n=0}^{\infty}M(n)$
is nonzero ($M(0)\ne 0$ by assumption) we see from the associativity
of vertex operators on $M$ that $\{u_n\sum_{n=0}^{\infty}M(n)|u\in
V^k,n\in \Z\}$ is a nonzero $V^0$-submodule of $M^k.$ Thus
$$\{u_n\sum_{n=0}^{\infty}M(n)|u\in V^k,n\in \Z\}=M^k.$$
In particular, $M^k$ is nonzero for all $k.$ 
Clearly $M^k$ and $M^i$ 
for $i\ne k$ are inequivalent $V^0$-modules. 

Let $0\ne w\in M^k.$ Then $\{u_nw|u\in V^i,n\in\Q\}$
is an admissible $V^0$-submodule of $M^{i+k}$ (where $i+k$ is understood modulo
$T$) and thus must be equal to $M^{i+k}$ for all $i.$ That
is, $\{u_nw|u\in V,n\in\Q\}=M$ and $M$ is an irreducible admissible 
$g$-twisted $V$-module.
\qed

As in the untwisted case, we can also formulate Proposition \ref{p6.3} 
in terms of theory of associative algebra $A_{g,n}(V)$ developed in [DLM3]
and [DLM5]. 

Let $V$ and $g$ be before. Fix $n=l+\frac{i}{T}\in\frac{1}{T}\Z$ with $l$ a nonnegative 
integer and $0\leq i\leq T-1.$ For $0\leq r\leq T-1$ we define 
$\delta_i(r)=1$ if $i\geq r$ and $\delta_i(r)=0$ if $i<r$.
We also set $\delta_i(T)=1.$ 
Let $O_{g,n}(V)$ be the linear span of all $u\circ_{g,n} v$ and $L(-1)u+L(0)u$
where for homogeneous $u\in V^r$ and $v\in V,$
$$u\circ_{g,n} v=\Res_{z}Y(u,z)v\frac{(1+z)^{\wt u-1+\delta_i(r)+l+r/T}}{z^{2l
+\delta_{i}(r)+\delta_{i}(T-r) }}.$$
Define the linear space $A_{g,n}(V)$ to be the quotient $V/O_{g,n}(V).$ Then
$A_{g,n}(V)$ is the untwisted associative algebra $A_n(V)$ as defined in
Section 4 if $g=1$ and is $A_g(V)$ in [DLM3] if $n=0.$ 
We also define a second product $*_{g,n}$ on $V$ for $u$ and $v$ as
above: 
$$u*_{g,n}v=\sum_{m=0}^{l}(-1)^m\binom{m+l}{l}\Res_zY(u,z)\frac{(1+z)^{\wt\,u+l}}{z^{l+m+1}}v$$
if $r=0$ and $u*_{g,n}v=0$ if $r>0.$
Extend linearly to obtain a bilinear product  on $V.$

Recall from [DLM3] that $V$ is called $g$-{\em rational} if any admissible  
$g$-twisted $V$-module is completely reducible. 
The following theorem was given in [DLM5].

\begin{thm} 
\label{theorem:2.2} Let $M=\sum_{n=0}^{\infty}M(n/T)$ be an admissible
$g$-twisted $V$-module with $M(0)\ne 0.$ We have

{\rm (i)} The product $*_{g,n}$ induces 
an associative algebra structure on $A_{g,n}(V)$ 
with the identity $\1+O_{g,n}(V).$ Moreover $\o+O_{g,n}(V)$ is a central
element of $A_{g,n}(V).$  

{\rm (ii)} The identity map on $V$ induces an onto algebra homomorphism
from $A_{g,n}(V)$ to $A_{g,m}(V)$ for $m,n\in\frac{1}{T}\Z$ and $0\leq m\leq n.$ 

{\rm (iii)} The map $u\mapsto o(u)$ gives a representation of $A_{g,n}(V)$ on $M(i)$ for $i\in \frac{1}{T}\Z$ and
$0\leq i\leq n.$ 
Moreover, $V$ is $g$-rational if and only if $A_{g,n}(V)$ are 
finite dimensional semisimple algebra for all $n.$ 
\end{thm}

Clearly, both the actions of $A_{g,n}(V)$ and $S_M(V^0)$ on 
$\sum_{0\leq m\leq n}M(m)$ are induced by $v\mapsto o(v).$ 
Combining Proposition \ref{p6.3} and Theorem \ref{theorem:2.2} gives 
an analogue of Theorem \ref{t6}.

\begin{thm}\label{t7} Assume that the Virasoro element $\omega$ of $V$ 
is nonzero. Then an admissible
 $g$-twisted $V$-module $M=\oplus_{n\geq 0}M(n)$ with $M(0)\ne 0$
is irreducible if and only if each $M(n)$ is an irreducible 
$A_{g,n}(V)$-module.
\end{thm}

\end{document}